\documentclass[11pt]{article}
\usepackage{amsmath,amssymb,amsfonts,latexsym}
\usepackage{amsthm}
\usepackage{mathrsfs} 
\usepackage{graphicx} 
\title{
Gradient-based iterative algorithms for solving Sylvester tensor equations
and the associated tensor nearness problems
}
\author{Mao-lin Liang}
\author{\small{
Maolin Liang\,\footnote{Corresponding author.
\ Email address:\,liangmaolin@tsnu.edu.cn/liangml2005@163.com\,(M. Liang),
bzheng@lzu.edu.cn (B. Zheng)
 }\ ,\,\,
Bing Zheng
}
\medskip\\
\small{$\textrm{School}$ of Mathematics and Statistics,
Lanzhou University,}\\
\small{Lanzhou, Gansu 730000, China}
}

\date{}

\begin{document}
\newtheorem{Proposition}{Proposition}[section]
\newtheorem{Lemma}{Lemma}[section]
\newtheorem{Theorem}{Theorem}[section]
\newtheorem{Example}{Example}[section] 
\newtheorem{Corollary}{Corollary}

\theoremstyle{definition}
\newtheorem{Definition}{Definition}[section]
\newtheorem{remark}{Remark}[section]
\maketitle

\begin{abstract}
   In this paper, an iterative algorithm is presented for solving
Sylvester tensor equation
$\mathscr{A}*_M\mathscr{X}+\mathscr{X}*_N\mathscr{C}=\mathscr{D}$,
where $\mathscr{A}$, $\mathscr{C}$ and $\mathscr{D}$ are given tensors with appropriate sizes,
and the symbol $*_N$ denotes the Einstein product.
By this algorithm, the solvability of this tensor equation can be determined
automatically, and the solution of which (when it is solvable) can be derived within
finite iteration steps for any initial iteration tensors in the absence of roundoff errors.
Particularly, the least F-norm solution of the aforementioned equation can be derived
by choosing special initial iteration tensors.
As application, we apply the proposed algorithm to the tensor nearness problem
related to the Sylvester tensor equation mentioned above.
It is proved that the solution to this problem can also be obtained within finite iteration steps
by solving another Sylvester tensor equation.
The performed numerical experiments show that the algorithm we propose here is promising.
\end{abstract}
\section{\bf Introduction}

   Tensors are multi-dimensional arrays \cite{Kolda}.
   An $Nth$-order and $I_1\times I_2 \times \dots \times I_N$-dimensional tensor
over the real field $\mathbb{R}$, consisting of
$I_1I_2 \cdots I_N$ entries, can be represented as
$$\mathscr{A}=(\mathscr{A}_{i_1\dots i_N})\
\textrm{with}\ \mathscr{A}_{i_1\dots i_N}\in \mathbb{R},\ 1\leq i_k \leq I_k,\ k=1,2,\ldots,N.$$
The set of this kind of tensors is denoted by $\mathbb{R}^{I_1\times I_2 \times \dots \times I_N}$.
   For
$\mathscr{A}\in\mathbb{R}^{I_1\times \dots \times I_M \times J_1\times \dots \times J_N}$
and
$\mathscr{B}\in\mathbb{R}^{J_1\times \dots \times J_N \times K_1\times \dots \times K_L}$,
the \textit{Einstein product} \cite{Einstein} of $\mathscr{A}$ and $\mathscr{B}$,
denoted by $\mathscr{A}*_N\mathscr{B}$,
is defined by the operation $*_N$ via
$$(\mathscr{A}*_N\mathscr{B})_{i_1 \dots i_M k_{1} \dots k_L}
=\sum_{j_1,\dots,j_N}a_{i_1 \dots i_M j_1 \dots j_N}
b_{j_1 \dots j_N k_1 \dots\ k_L}.$$
   Tensor models are employed in numerous disciplines addressing the problem of finding
multilinear structure in multiway data-sets. In particular, tensor equations with Einstein product
model many phenomena in engineering and science, including continuum physics and engineering, isotropic
and anisotropic elasticity \cite{Li1,Li2,Lai}.
   For example, by using the central difference approximation, the three-dimensional Poisson equations
can be discretized as the following multilinear system \cite{Brazell}
$$\mathscr{A}*_3\mathscr{X}=\mathscr{B},\ \mathscr{X} \in\mathbb{R}^{N\times N \times N},$$
where tensors
$\mathscr{A}\in\mathbb{R}^{N\times N \times N \times N\times N \times N}$,
$\mathscr{B}\in\mathbb{R}^{N\times N \times N}$.
   The general form of the above tensor equation is as follows:
\begin{equation}\label{eqA01}
\mathscr{A}*_M\mathscr{X}=\mathscr{B},\\
\end{equation}
where
$\mathscr{A}\in\mathbb{R}^{K_1\times \dots \times K_P \times I_1\times \dots \times I_M}$ and
$\mathscr{B}\in\mathbb{R}^{K_1\times \dots \times K_P \times J_1\times \dots \times J_N}$ are given tensors,
and
$\mathscr{X}\in\mathbb{R}^{I_1\times \dots \times I_M\times J_1\times \dots \times J_N}$ is unknown.
   Brazell et. al \cite{Brazell} researched the tensor equation \eqref{eqA01} and the associated least-square
problem by introducing the notion of
inverse or pseudo-inverse of a tensor.
   Recently, Sun et. al. \cite{Sun} extended the inverse in \cite{Brazell} and put forward
the concept of Moore-Penrose inverses of tensors which provides the way
to represent the general solution of the tensor equation \eqref{eqA01}
in the sense that it is consistent (namely, there exists a tensor $\mathscr{X}^*$ satisfying \eqref{eqA01}).
   Besides, the authors also considered the Sylvester tensor equation
\begin{equation}\label{eqA02}
\mathscr{A}*_M\mathscr{X}+\mathscr{X}*_N\mathscr{C}=\mathscr{D},\\
\end{equation}
in which
$\mathscr{A} \in\mathbb{R}^{I_1\times \cdots \times I_M \times I_1\times \cdots \times I_M}$,
$\mathscr{C} \in\mathbb{R}^{J_1\times \cdots \times J_N \times J_1\times \cdots \times J_N}$ and
$\mathscr{D} \in\mathbb{R}^{I_1\times \cdots \times I_M \times J_1\times \cdots \times J_N}$
are given tensors, and
$\mathscr{X} \in \mathbb{R}^{I_1\times \dots \times I_M \times J_1\times \dots \times J_N}$
is the one to be determined.

   This equation is a generalization of the well-known Sylvester matrix equation,
comes from the finite element, finite difference or spectral method \cite{Li1,Li2},
and plays an important role in discretization of linear partial differential equations in high dimension \cite{Lai,Itskov,Brazell,Sun}.
   Based on the operations of `block tensors',
it is proved \cite{Sun} that \eqref{eqA02} is equivalent to
\begin{equation*}
[\mathscr{A}\ \mathscr{I}_1]*_M
\left[
  \begin{array}{cc}
   \mathscr{X} & \mathscr{O}\\
   \mathscr{O} & \mathscr{X}\\
  \end{array}
\right]*N
\left[
  \begin{array}{c}
   \mathscr{I}_2\\
   \mathscr{C}\\
  \end{array}
\right]
=\mathscr{D},\\
\end{equation*}
where $\mathscr{I}_i$ $(i=1,2)$ are the identity tensors with appropriate size.
Nevertheless, it is difficult to derive the explicit expression of the solution
via the Moore-Penrose inverse for the last tensor equation,
since it requires the $2\times 2$ block structure.
   On the other hand, one could compute the exact solution of such an equation
by converting it into the form of \eqref{eqA01}
via the Kronecker product,
but the computational efforts rapidly increases with the dimensions of the tensors.

   The purpose of this paper is to solve the Sylvester tensor equation \eqref{eqA02}
by establishing the gradient-based iterative method
twisted from the ones given in \cite{Brazell,Wang,Peng,Chenz,Liang},
see Section 3 for details.
%
%
%
   It is theoretically shown that the proposed approach can be capable of finding
the solution of \eqref{eqA02}
within finite iteration steps  for any initial iteration tensors.
Especially, the least Frobenius norm (F-norm for short) solution of which can also be derived
by choosing appropriate initial iteration tensors.

   Another problem we are interested in is the following constrained minimization problem
related to the Sylvester tensor equation \eqref{eqA02}:
\begin{equation}\label{eqA03}
\min\limits_{\mathscr{A}*_M\mathscr{X}+\mathscr{X}*_N\mathscr{C}=\mathscr{D}}
\|\mathscr{X}-\mathscr{X}_0\|,\\
\end{equation}
where $\mathscr{X}_0\in \mathbb{R}^{I_1\times \dots \times I_M\times J_1\times \dots \times J_N}$
is a given tensor, the symbol $\|\cdot\|$ denotes the Frobenius norm of a tensor.
    This problem is a natural generalization of the matrix nearness problem
\cite{Higham,Liao,Dhillon,Noschese},
low rank approximation problem \cite{Zhang,Silva,Qi} and
tensor completion problem \cite{Liu,Gandy,Geng} equipped with F-norm and multilinear constraints.
We call \eqref{eqA03} the tensor nearness problem.
   Under certain conditions,
it will be proved that the solution to the tensor nearness problem \eqref{eqA03}
is unique, and can be gained by applying the proposed algorithm to
another Sylvester tensor equation, see Section 4 for details.
   Particularly, when the tensor $\mathscr{C}$ in \eqref{eqA03} vanishes,
we have proved that the unique solution to the tensor nearness problem
can be represented by means of the Moore-Penrose inverses of the known tensors.
Nevertheless, it is well-known that it is not easy to find the Moore-Penrose inverse of a tensor.
Therefore, the work of this paper avoids this curse
and thus can be regarded as a continuation of \cite{Liang2018}.

   The remainder of this paper is organized as follows.
   Section 2 reviews some notation and definitions related to tensors.
   Section 3 contains the gradient-based iterative algorithm for solving the tensor equation \eqref{eqA02} as well as its convergence analysis.
   Section 4 is devoted to addressing the tensor nearness problem \eqref{eqA03}.
   Section 5 provides some numerical examples to illustrate the efficiency
of the proposed iterative algorithms.
   Finally, a conclusion is appended to end this paper.

\section{Preliminaries}

    Throughout this paper,
tensors are denoted by calligraphic letters, e.g., $\mathscr{A}$, $\mathscr{B}$, $\mathscr{C}$;
matrices are denoted by boldface capital letters, e.g., $\textbf{A}$, $\textbf{B}$, $\textbf{C}$;
Vectors are denoted by boldface lowercase letters, e.g., $\textbf{a}$, $\textbf{b}$, $\textbf{c}$;
Scalars are denoted by lowercase letters, e.g., $a$, $b$, $c$.
   For a higher-order tensor, subtensors are formed when a subset of the indices is fixed, and
a colon is used to indicate all elements of a mode.
   For example, if a tensor $\mathscr{A}\in\mathbb{R}^{I\times J\times K}$,
its column, row, and tube fibers, which are denoted
by $\mathscr{A}(:,j,k),\ \mathscr{A}(i,:,k)$ and $\mathscr{A}(i,j,:)$, respectively.
   Moreover, the horizontal, lateral, and frontal slices are represented by
$\mathscr{A}(i,:,:),\ \mathscr{A}(:,j,:)$ and $\mathscr{A}(:,:,k)$, respectively.

   The following definitions and conclusions will be used later.
\begin{Definition}\label{def01} (\cite{Brazell})
   For
$\mathscr{A}=(\mathscr{A}_{i_1 \ldots i_M j_1 \ldots j_N})
\in\mathbb{R}^{I_1\times \dots \times I_M \times J_1\times \dots \times J_N}$,
its transpose, denoted by $\mathscr{A}^T$,
is a $J_1\times \dots \times J_N\times I_1\times \dots \times I_M$ tensor with the entries
$\widehat{\mathscr{A}}_{i_1 \ldots i_N j_1 \ldots j_M}
=\mathscr{A}_{j_1 \ldots j_M i_1 \ldots i_N }$.

   Particularly, if $\mathscr{A}=(\mathscr{A}_{i_1 \ldots i_M j_1 \ldots j_M})
\in\mathbb{R}^{I_1\times \dots \times I_M \times I_1\times \dots \times I_M}$,
the trace of $\mathscr{A}$, denoted by $tr(\mathscr{A})$, is defined as
$tr(\mathscr{A})=\sum\limits_{i_1,\ldots, i_M}\mathscr{A}_{i_1 \ldots i_M i_1 \ldots i_M}$.
\end{Definition}

   By Definition \ref{def01}, the inner product of two tensors
$\mathscr{A},\mathscr{B}
\in\mathbb{R}^{I_1\times \dots \times I_M \times J_1\times \dots \times J_N}$
is defined by
$<\mathscr{A},\ \mathscr{B}>=tr(\mathscr{B}^T*_N \mathscr{A})$,
which induces the Frobenius norm of a tensor, i.e.,
$\|\mathscr{A}\|=\sqrt{<\mathscr{A},\ \mathscr{A}>}$.
   Especially, if $<\mathscr{A},\ \mathscr{B}>=0$, we say that the two tensors are orthogonal each other.
   Moreover, it is easy to verify that the following results hold true.
\begin{Lemma}\label{le2.01}
   Let $\mathscr{A}, \mathscr{B}, \mathscr{C}
\in\mathbb{R}^{I_1\times \cdots \times I_M \times I_1\times \cdots \times I_M}$
and $\alpha,\beta \in \mathbb{R}$, then
\item \rm{(I)} $tr(\alpha\cdot \mathscr{A}+\beta\cdot \mathscr{B})
=\alpha\cdot tr(\mathscr{A})+\beta\cdot tr(\mathscr{B})$;
\item \rm{(II)} $tr(\mathscr{A}*_M \mathscr{B}*_M \mathscr{C})
=tr(\mathscr{B}*_M \mathscr{C}*_M \mathscr{A})
=tr(\mathscr{C}*_M \mathscr{A}*_M \mathscr{B})$.
\end{Lemma}

\begin{Definition}\label{def02} 
   For a tensor
$\mathscr{A} \in\mathbb{R}^{I_1\times \cdots \times I_M \times J_1\times \cdots \times J_N}$,
$\textrm{Vec}(\mathscr{A})\in\mathbb{R}^{(I_1\cdot \ldots \cdot I_M) \times J_1\times \cdots \times J_N}$
is obtained by lining up all the subtensors, $\mathscr{A}(i_1,\ldots,i_M,:,\ldots,:)$
with $1\leq i_j\leq I_j$ and $j=1,2,\ldots,M$, in a column;
e.g., the $k$th subblock of $\mathscr{A}$ is the subtensor $\mathscr{A}(i_1,\ldots,i_M,:,\ldots,:)$
satisfying $k=\rm{ivec}(\emph{\rm{\textbf{i}}},\mathbb{I})$, where
$\rm{ivec}(\cdot)$ is the index mapping function \cite{Ragnarsson}, i.e.,\\
\centerline{
$\rm{ivec}(\emph{\rm{\textbf{i}}},\mathbb{I})$
$:=i_1+\sum\limits_{j=2}^{M} (i_j-1) \prod\limits_{s=1}^{j-1} I_s$ and $\mathbb{I}=\{I_1,\ldots,I_M\}$.}
\end{Definition}

   Specifically, if $\mathscr{A} \in\mathbb{R}^{2\times 2\times 2 \times J_1\times \dots \times J_N}$,
then
\begin{equation*}
\textrm{Vec}(\mathscr{A})=
\left[
  \begin{array}{c}
  \mathscr{A}(1,1,1,:,\ldots,:)\\
  \mathscr{A}(2,1,1,:,\ldots,:)\\
  \mathscr{A}(1,2,1,:,\ldots,:)\\
  \mathscr{A}(2,2,1,:,\ldots,:)\\
  \mathscr{A}(1,1,2,:,\ldots,:)\\
  \mathscr{A}(2,1,2,:,\ldots,:)\\
  \mathscr{A}(1,2,2,:,\ldots,:)\\
  \mathscr{A}(2,2,2,:,\ldots,:)\\
  \end{array}
\right].
\end{equation*}
    We should mention that the definition of Vec is slightly different from that given in \cite{Sun}.

\begin{Definition}\label{def03} (\cite{Sun})
   The Kronecker product of
$\mathscr{A}=(\mathscr{A}_{i_1 \ldots i_M j_1 \ldots j_N})
\in\mathbb{R}^{I_1\times \cdots \times I_M \times J_1\times \cdots \times J_N}$ and
$\mathscr{B} \in\mathbb{R}^{K_1\times \cdots \times K_P \times L_1\times \cdots \times L_Q}$,
denoted by
$\mathscr{A}\otimes \mathscr{B}$,
is a `Kr-block tensor', whose ($r,s$)-subblock is
$(\mathscr{A}_{i_1 \ldots i_M j_1 \ldots j_N}\mathscr{B})$ in which
$r=\rm{ivec}(\emph{\rm{\textbf{i}}},\mathbb{I})$
and $s=\rm{ivec}(\emph{\rm{\textbf{j}}},\mathbb{J})$ for $\mathbb{J}=\{J_1,\ldots,J_N\}$.
\end{Definition}
   The Kronecker product of tensors has the following basic properties:
\begin{Lemma}\label{le2.1}\rm{(\cite{Sun,Behera})}
   Let
$\mathscr{A} \in\mathbb{R}^{I_1\times \cdots \times I_M \times I_1\times \cdots \times I_M}$,
$\mathscr{B} \in\mathbb{R}^{I_1\times \cdots \times I_M \times J_1\times \cdots \times J_N}$,\\
$\mathscr{C} \in\mathbb{R}^{J_1\times \cdots \times J_N \times J_1\times \cdots \times J_N}$ and
$\mathscr{D} \in\mathbb{R}^{I_1\times \cdots \times I_M \times J_1\times \cdots \times J_N}$. Then
\item \rm{(I)} $(\mathscr{B}+\mathscr{D})^T=\mathscr{B}^T+\mathscr{D}^T$.
\item \rm{(II)} $(\mathscr{A}\otimes\mathscr{B})\otimes\mathscr{C}
=\mathscr{A}\otimes(\mathscr{B}\otimes\mathscr{C})$.
\item \rm{(III)} $(\mathscr{A}\otimes\mathscr{B})*_N(\mathscr{D}\otimes\mathscr{C})
=(\mathscr{A}*_M\mathscr{D})\otimes(\mathscr{B}*_N\mathscr{C})$.
\item \rm{(IV)} $\textrm{Vec}(\mathscr{A}*_M\mathscr{B}*_N\mathscr{C})
    =(\mathscr{C}^T\otimes\mathscr{A})*_N \textrm{Vec}(\mathscr{B})$.
\end{Lemma}
\begin{Definition}\label{def04}\rm{(\cite{Brazell})}
   Define the transformation $\psi$ from the tensor space
$\mathbb{R}^{I_1\times \cdots \times I_M \times J_1\times \cdots \times J_N}$
to the matrix space
$\mathbb{R}^{(I_1\cdot \ldots\cdot I_M) \times (J_1\cdot \ldots \cdot J_N)}$
as
\begin{equation*}
\begin{split}
\psi:\ \mathbb{C}^{I_1\times \cdots \times I_M \times J_1\times \cdots \times J_N}
    &\longrightarrow
    \mathbb{C}^{(I_1\cdot \ldots\cdot I_M) \times (J_1\cdot \ldots\cdot J_N)}\\
\mathscr{A}_{i_1\ldots i_M j_1\ldots j_N}
   &\longrightarrow
\textbf{A}_{\rm{ivec}(\emph{\rm{\textbf{i}}},\mathbb{I}),\ \rm{ivec}(\emph{\rm{\textbf{j}}},\mathbb{J})}.
\end{split}
\end{equation*}
\end{Definition}

   Obviously, the transformation $\psi$ is a bijection, which provides a way to unfold one tensor.
   For example, if
$\mathscr{A} \in\mathbb{R}^{3\times 3\times 3\times 3}$, each frontal slice
$\mathscr{A}(:,:,k,l)$ with $k,l=1,2,3$ is a $3\times 3$ matrix.
If partition the modes of the tensor $\mathscr{A}$ from the middle, then the vector
$\textrm{vec}(\mathscr{A}(:,:,k,l))$ corresponds to the $[k+3(l-1)]$th column of the unfolding matrix
$\textbf{A}=\psi(\mathscr{A})$, that is,
$$
\textbf{A}=
\left[
   \begin{array}{ccccccccc}
\mathscr{A}_{1111}\vline& \mathscr{A}_{1121}\vline&\mathscr{A}_{1131}\vline&\mathscr{A}_{1112}\vline&\mathscr{A}_{1122}\vline
&\mathscr{A}_{1132}\vline&\mathscr{A}_{1113}\vline&\mathscr{A}_{1123}\vline&\mathscr{A}_{1133}\\
\mathscr{A}_{2111}\vline& \mathscr{A}_{2121}\vline&\mathscr{A}_{2131}\vline&\mathscr{A}_{2112}\vline&\mathscr{A}_{2122}\vline
&\mathscr{A}_{2132}\vline&\mathscr{A}_{2113}\vline&\mathscr{A}_{2123}\vline&\mathscr{A}_{2133}\\
\mathscr{A}_{3111}\vline& \mathscr{A}_{3121}\vline&\mathscr{A}_{3131}\vline&\mathscr{A}_{3112}\vline&\mathscr{A}_{3122}\vline
&\mathscr{A}_{3132}\vline&\mathscr{A}_{3113}\vline&\mathscr{A}_{3123}\vline&\mathscr{A}_{3133}\\
\mathscr{A}_{1211}\vline& \mathscr{A}_{1221}\vline&\mathscr{A}_{1231}\vline&\mathscr{A}_{1212}\vline&\mathscr{A}_{1222}\vline
&\mathscr{A}_{1232}\vline&\mathscr{A}_{1213}\vline&\mathscr{A}_{1223}\vline&\mathscr{A}_{1233}\\
\mathscr{A}_{2211}\vline& \mathscr{A}_{2221}\vline&\mathscr{A}_{2231}\vline&\mathscr{A}_{2212} \vline&\mathscr{A}_{2222}\vline&\mathscr{A}_{2232}\vline&\mathscr{A}_{2213}\vline
&\mathscr{A}_{2223}\vline&\mathscr{A}_{2233}\\
\mathscr{A}_{3211}\vline& \mathscr{A}_{3221} \vline& \mathscr{A}_{3231}\vline&\mathscr{A}_{3212}\vline&\mathscr{A}_{3222}\vline
&\mathscr{A}_{3232}\vline&\mathscr{A}_{3213}\vline&\mathscr{A}_{3223}\vline&\mathscr{A}_{3233} \\
\mathscr{A}_{1311}\vline&  \mathscr{A}_{1321}\vline& \mathscr{A}_{1331}\vline&\mathscr{A}_{1312}\vline&\mathscr{A}_{1322}\vline
&\mathscr{A}_{1332}\vline&\mathscr{A}_{1313}\vline&\mathscr{A}_{1323}\vline&\mathscr{A}_{1333} \\
\mathscr{A}_{2311}\vline&  \mathscr{A}_{2321}\vline& \mathscr{A}_{2331}\vline
&\mathscr{A}_{2312}\vline& \mathscr{A}_{2322}\vline& \mathscr{A}_{2332}\vline
&\mathscr{A}_{2313}  \vline& \mathscr{A}_{2323}\vline&\mathscr{A}_{2333}\\
\mathscr{A}_{3311}\vline&  \mathscr{A}_{3321}\vline& \mathscr{A}_{3331}\vline
&\mathscr{A}_{3312}\vline& \mathscr{A}_{3322}\vline& \mathscr{A}_{3332}\vline
&\mathscr{A}_{3313}\vline&\mathscr{A}_{3323}\vline&\mathscr{A}_{3333}\\
   \end{array}
 \right].
$$

   From the definition of $\psi$, one can observe
that the entry $\mathscr{A}_{i_1\dots i_M j_1 \dots j_N}$ of the tensor
$\mathscr{A}\in\mathbb{R}^{I_1\times \cdots \times I_M \times J_1\times \cdots\times J_N}$ is exactly
the ($\rm{ivec}(\textrm{\textbf{i}}, \mathbb{I}),\rm{ivec}(\textrm{\textbf{j}}, \mathbb{J})$)-element
of the image matrix $\psi(\mathscr{A})$.
Thus, the \textit{identity} \textit{tensor} of size
$I_1\times \cdots \times I_M \times I_1\times \cdots \times I_M$,
denoted by $\mathscr{I}$, consists of the entries
\begin{equation*}
\begin{aligned}
\mathscr{I}_{i_1\dots i_Mj_1\dots j_M}=\prod_{k=1}^M\delta_{i_kj_k} \
\textrm{with}\ \delta_{i_k j_k}=\left\{
                       \begin{array}{c}
                         1 ,\quad \textrm{if} \ i_k=j_k,\\
                         0 ,\quad \textrm{if} \ i_k\neq j_k.\\
                       \end{array}
                     \right.
\end{aligned}
\end{equation*}

   In addition, for tensor
$\mathscr{A}\in\mathbb{R}^{I_1\times \dots \times I_M \times J_1\times \dots \times J_N}$,
its range space is defined by
$$R(\mathscr{A})=\{\ \mathscr{Y}\ | \ \mathscr{Y}=\mathscr{A}*_N \mathscr{X},\
\forall\ \mathscr{X}\in
\mathbb{R}^{J_1\times \dots \times J_N} \}.$$

\section{The iterative algorithm and its convergence analysis}

   In this section, we propose the gradient-based iterative algorithm for solving
the Sylvester tensor equation \eqref{eqA02}, and then analyze its convergence.

   The iterative algorithm for solving \eqref{eqA02} is described as below:

\noindent\textbf{Algorithm 3.1}
\medskip\\
\noindent $\textsf{Step\,\,\,1}$: Input
$\mathscr{A} \in\mathbb{R}^{I_1\times \cdots \times I_M \times I_1\times \cdots \times I_M}$,
$\mathscr{C} \in\mathbb{R}^{J_1\times \cdots \times J_N \times J_1\times \cdots \times J_N}$,
$\mathscr{D} \in\mathbb{R}^{I_1\times \cdots \times I_M \times J_1\times \cdots \times J_N}$,
\\
\indent\indent\,
and an initial iteration tensor
$\mathscr{X}^{(1)}\in\mathbb{R}^{I_1\times \dots \times I_M \times J_1\times \dots \times J_N}$.

\noindent $\textsf{Step\,\,\,2}$: Compute
$\mathscr{R}^{(1)}=\mathscr{D}-\mathscr{A}*_M \mathscr{X}^{(1)}-\mathscr{X}^{(1)}*_N \mathscr{C}$, and
\\
\indent\indent\,
$\mathscr{P}^{(1)}=\mathscr{A}^T *_M \mathscr{R}^{(1)}+\mathscr{R}^{(1)}*_N \mathscr{C}^T$.

\noindent $\textsf{Step\,\,\,3}$: Compute
%
$\mathscr{X}^{(k+1)}=\mathscr{X}^{(k)}
+\dfrac{\|\mathscr{R}^{(k)}\|^2}{\|\mathscr{P}^{(k)}\|^2}\ \mathscr{P}^{(k)}$.
\medskip\\
\noindent
$\textsf{Step\,\,\,4}$: Compute
$\mathscr{R}^{(k+1)}=\mathscr{D}-\mathscr{A}*_M \mathscr{X}^{(k+1)}-\mathscr{X}^{(k+1)}*_N \mathscr{C}$,
 and\\
%
%
\indent\indent\indent\indent\indent
$\mathscr{P}^{(k+1)}=\mathscr{A}^T*_M \mathscr{R}^{(k+1)}+\mathscr{R}^{(k+1)}*_N \mathscr{C}^T
+\dfrac{\|\mathscr{R}^{(k+1)}\|^2}{\|\mathscr{R}^{(k)}\|^2}\ \mathscr{P}^{(k)}$.

\indent\indent\,
If $\mathscr{R}^{(k+1)}=0$, or $\mathscr{R}^{(k+1)}\neq 0$, $\mathscr{P}^{(k)}=0$, stop;
Otherwise, $k:=k+1$, goto Step 3.

  In what follows, we show that the sequence $\{\mathscr{X}^{(k)}\}$ generated by Algorithm 3.1
converges to a solution of \eqref{eqA02} within finite iteration steps in the absence of roundoff errors
for any initial iteration tensor $\mathscr{X}^{(1)}$.
  For ease of expression, denote
  $$\alpha(k):=\dfrac{\|\mathscr{R}^{(k)}\|^2}{\|\mathscr{P}^{(k)}\|^2},\ \
  \beta(k+1):=\dfrac{\|\mathscr{R}^{(k+1)}\|^2}{\|\mathscr{R}^{(k)}\|^2}.$$

\begin{Lemma}\label{le3.1} Let $\{\mathscr{R}^{(i)}\}$ and $\{\mathscr{P}^{(i)}\}$ $(i=1,2,\ldots\,)$
be the sequences generated by Algorithm 3.1, then, for $j\geq 2$, it holds that
\begin{equation}\label{eqC01}
\begin{split}
tr\left({\mathscr{R}^{(i+1)}}^T *_M \mathscr{R}^{(j)}\right)
=&tr\left({\mathscr{R}^{(i)}}^T *_M \mathscr{R}^{(j)}\right)
-\alpha(i) \cdot tr\left({\mathscr{P}^{(i)}}^T *_M \mathscr{P}^{(j)}\right)\\
&+\,\alpha(i) \cdot \beta(j) \cdot
      tr\left({\mathscr{P}^{(i)}}^T *_M \mathscr{P}^{(j-1)}\right).\\
\end{split}
\end{equation}
\end{Lemma}

\noindent\textit{Proof}. By Algorithm 3.1, we have
\begin{equation*}\label{eqC03}
\begin{split}
&tr\left({\mathscr{R}^{(i+1)}}^T *_M \mathscr{R}^{(j)}\right)\\
&=tr\left((\mathscr{R}^{(i)}-\alpha(i) \cdot
(\mathscr{A}*_M \mathscr{P}^{(k)}+\mathscr{P}^{(i)}*_N \mathscr{C}))^T*_M \mathscr{R}^{(j)}\right)\\
&=tr\left({\mathscr{R}^{(i)}}^T *_M \mathscr{R}^{(j)}\right) -\alpha(i) \cdot
tr\left((\mathscr{A}*_M \mathscr{P}^{(i)}+\mathscr{P}^{(i)}*_N \mathscr{C})^T *_M \mathscr{R}^{(j)}\right)\\
&=tr\left({\mathscr{R}^{(i)}}^T *_M \mathscr{R}^{(j)}\right) -\alpha(i) \cdot
tr\left({\mathscr{P}^{(i)}}^T *_M (\mathscr{A}^T *_M \mathscr{R}^{(j)}+\mathscr{R}^{(j)}*_N \mathscr{C}^T)\right)\\
&=tr\left({\mathscr{R}^{(i)}}^T *_M \mathscr{R}^{(j)}\right) -\alpha(i) \cdot
tr\left({\mathscr{P}^{(i)}}^T *_M (\mathscr{P}^{(j)} -\beta(j) \cdot \mathscr{P}^{(j-1)})\right),
\end{split}
\end{equation*}
which implies that the equality \eqref{eqC01} holds true. $\square$

\indent The next lemma reveals the orthogonality of the sequences
$\{\mathscr{R}^{(i)}\}$ and $\{\mathscr{P}^{(i)}\}$ generated by Algorithm 3.1,
which is similar to the classical conjugate gradient method \cite{Golub}.

\begin{Lemma}\label{le3.2} Let $\{\mathscr{R}^{(i)}\}$ and $\{\mathscr{P}^{(i)}\}$ $(i=1,2,\ldots\,)$
be the sequences generated by Algorithm 3.1. Then
\begin{equation}\label{eqC04}
tr\left({\mathscr{R}^{(i)}}^T *_M \mathscr{R}^{(j)}\right)=0, \ \
tr\left({\mathscr{P}^{(i)}}^T *_M \mathscr{P}^{(j)}\right)=0, \ \
i,j=1,2,\ldots,t\, (t\geq 2), \, i\neq j.\\
\end{equation}
\end{Lemma}

\noindent\textit{Proof}.  We prove \eqref{eqC04} by induction.
Since
$tr\left({\mathscr{R}^{(i)}}^T *_M \mathscr{R}^{(j)}\right)
=tr\left({\mathscr{R}^{(j)}}^T *_M \mathscr{R}^{(i)}\right)$, so we only consider
the case: $i\geq j$.
\medskip \\
\indent When $t=2$, from Algorithm 3.1, we obtain
\begin{equation}\label{eqC05}
\begin{split}
&tr\left({\mathscr{R}^{(2)}}^T *_M \mathscr{R}^{(1)}\right)\\
&=tr\left({\mathscr{R}^{(1)}}^T *_M \mathscr{R}^{(1)}\right)
-\alpha(1)\cdot tr\left((\mathscr{A}*_M \mathscr{P}^{(1)}+\mathscr{P}^{(1)}*_N \mathscr{C})^T *_M \mathscr{R}^{(1)}\right)\\
&=tr\left({\mathscr{R}^{(1)}}^T *_M \mathscr{R}^{(1)}\right)
-\alpha(1)
\cdot tr\left({\mathscr{P}^{(1)}}^T *_M (\mathscr{A}^T *_M \mathscr{R}^{(1)}+\mathscr{R}^{(1)}*_N \mathscr{C}^T)\right)\\
&=tr\left({\mathscr{R}^{(1)}}^T *_M \mathscr{R}^{(1)}\right)
-\alpha(1)
\cdot tr\left({\mathscr{P}^{(1)}}^T *_M \mathscr{P}^{(1)}\right)\\
&=0,\\
\end{split}
\end{equation}
and
\begin{equation}\label{eqC06}
\begin{split}
&tr\left({\mathscr{P}^{(2)}}^T *_M \mathscr{P}^{(1)}\right)\\
&=tr\left((\mathscr{A}^T*_M \mathscr{R}^{(2)}+\mathscr{R}^{(2)}*_N \mathscr{C}^T
+\beta(2) \cdot \mathscr{P}^{(1)})^T *_M {\mathscr{P}}^{(1)}\right)\\
&=tr\left({\mathscr{R}^{(2)}}^T *_M (\mathscr{A}*_M \mathscr{P}^{(1)}+\mathscr{P}^{(1)}*_N \mathscr{C})\right)
+ \beta(2) \cdot tr\left({\mathscr{P}^{(1)}}^T *_M \mathscr{P}^{(1)}\right)\\
&=tr\left({\mathscr{R}^{(2)}}^T *_M (\mathscr{R}^{(1)}-\mathscr{R}^{(2)})
\cdot \dfrac{1}{\alpha(1)}\right)
+ \beta(2) \cdot tr\left({\mathscr{P}^{(1)}}^T *_M \mathscr{P}^{(1)}\right)\\
&=-\dfrac{1}{\alpha(1)} \cdot tr\left({\mathscr{R}^{(2)}}^T *_M \mathscr{R}^{(2)}\right)
+ \beta(2) \cdot tr\left({\mathscr{P}^{(1)}}^T *_M \mathscr{P}^{(1)}\right)\\
&=0.\\
\end{split}
\end{equation}

   Suppose that \eqref{eqC04} holds for $t=s$, that is,
$$tr\left({\mathscr{R}^{(s)}}^T *_M \mathscr{R}^{(j)}\right)=0,\,
tr\left({\mathscr{P}^{(s)}}^T *_M \mathscr{P}^{(j)}\right)=0,\, j=1,2,\dots,s-1.$$
In view of Lemma \ref{le3.1}, when $t=s+1$, we have
\begin{equation}\label{eqC07}
\begin{split}
tr\left({\mathscr{R}^{(s+1)}}^T *_M \mathscr{R}^{(s)}\right)
&=tr\left({\mathscr{R}^{(s)}}^T *_M \mathscr{R}^{(s)}\right)
-\,\alpha(s) \cdot tr\left({\mathscr{P}^{(s)}}^T *_M \mathscr{P}^{(s)}\right)\\
&\ \ \ +\alpha(s) \cdot \beta(s) \cdot tr\left({\mathscr{P}^{(s)}}^T *_M \mathscr{P}^{(s-1)}\right)\\
&=tr\left({\mathscr{R}^{(s)}}^T *_M \mathscr{R}^{(s)}\right)
-\alpha(s) \cdot tr\left({\mathscr{P}^{(s)}}^T *_M \mathscr{P}^{(s)}\right)\\
&=0,\\
\end{split}
\end{equation}
and
\begin{equation}\label{eqC08}
\begin{split}
&tr\left({\mathscr{P}^{(s+1)}}^T *_M \mathscr{P}^{(s)}\right)\\
&=tr\left((\mathscr{A}^T*_M \mathscr{R}^{(s+1)}+\mathscr{R}^{(s+1)}*_N \mathscr{C}^T
   +\beta(s+1) \cdot \mathscr{P}^{(s)})^T *_M {\mathscr{P}}^{(s)}\right)\\
&=tr\left({\mathscr{R}^{(s+1)}}^T *_M (\mathscr{R}^{(s)}-\mathscr{R}^{(s+1)})
\cdot \dfrac{1}{\alpha(s)}\right)+ \beta(s+1) \cdot tr\left({\mathscr{P}^{(s)}}^T *_M \mathscr{P}^{(s)}\right)\\
&=-\dfrac{1}{\alpha(s)} \cdot tr\left({\mathscr{R}^{(s+1)}}^T *_M \mathscr{R}^{(s+1)}\right)
+ \beta(s+1) \cdot tr\left({\mathscr{P}^{(s)}}^T *_M \mathscr{P}^{(s)}\right)\\
&=0.\\
\end{split}
\end{equation}

    Now we consider the cases $j=1,2,\ldots,s-1$.
    In fact, when $j=1$, similar to the proofs of \eqref{eqC05} and \eqref{eqC06}, we have
\begin{equation}\label{eqC09}
\begin{split}
tr\left({\mathscr{R}^{(s+1)}}^T *_M \mathscr{R}^{(1)}\right)
&=tr\left((\mathscr{R}^{(s)}-\alpha(s) \cdot
                 (\mathscr{A}*_M \mathscr{P}^{(s)}+\mathscr{P}^{(s)}*_N \mathscr{C}))^T *_M \mathscr{R}^{(1)}\right)\\
&=-\,\alpha(s)
\cdot tr\left({\mathscr{P}^{(s)}}^T *_M (\mathscr{A}^T *_M \mathscr{R}^{(1)}
      +\mathscr{R}^{(1)} *_M \mathscr{C}^T)\right)\\
&=-\,\alpha(s)
\cdot tr\left({\mathscr{P}^{(s)}}^T *_M \mathscr{P}^{(1)}\right)\\
&=0,
\end{split}
\end{equation}
and
\begin{equation}\label{eqC10}
\begin{split}
&tr\left({\mathscr{P}^{(s+1)}}^T *_M \mathscr{P}^{(1)}\right)\\
&=tr\left((\mathscr{A}^T*_M \mathscr{R}^{(s+1)}+\mathscr{R}^{(s+1)}*_N \mathscr{C}^T
   +\beta(s+1) \cdot \mathscr{P}^{(s)})^T *_M {\mathscr{P}}^{(1)}\right)\\
&=tr\left({\mathscr{R}^{(s+1)}}^T *_M (\mathscr{A}*_M \mathscr{P}^{(1)}+\mathscr{P}^{(1)}*_N \mathscr{C})\right)\\
&=-\dfrac{1}{\alpha(1)} \cdot
tr\left({\mathscr{R}^{(s+1)}}^T *_M (\mathscr{R}^{(1)}-\mathscr{R}^{(2)})\right)\\
&=0.
\end{split}
\end{equation}

\noindent When $2 \leq j \leq s-1$, similar to the proofs of \eqref{eqC07} and \eqref{eqC08},
using Lemma 3.1 once again, we can respectively deduce that
\begin{equation*}
tr\left({\mathscr{R}^{(s+1)}}^T *_M \mathscr{R}^{(j)}\right)=0 \ \textrm{and}\
tr\left({\mathscr{P}^{(s+1)}}^T *_M \mathscr{P}^{(j)}\right)=0,\\
\end{equation*}
which, together with \eqref{eqC05}-\eqref{eqC10}, indicates that \eqref{eqC04} holds. $\square$
\begin{Lemma}\label{le3.3} Suppose that $\widetilde{\mathscr{X}}$ is an arbitrary solution
of the tensor equation \eqref{eqA02}, then the sequences $\{\mathscr{R}^{(k)}\}$
and $\{\mathscr{P}^{(k)}\}$ satisfy
\begin{equation}\label{eqC11}
tr\left((\widetilde{\mathscr{X}}-\mathscr{X}^{(k)})^T *_M \mathscr{P}^{(k)}\right)
=\|\mathscr{R}^{(k)}\|^2,\ \ k=1,2,\ldots.\\
\end{equation}
\end{Lemma}
\noindent\textit{Proof}. We prove \eqref{eqC11} by induction as well.
   When $k=1$, it follows from Algorithm 3.1 and Lemma 3.2 that
\begin{equation}\label{eqC12}
\begin{split}
&tr\left((\widetilde{\mathscr{X}}-\mathscr{X}^{(1)})^T *_M \mathscr{P}^{(1)}\right)\\
&=tr\left((\widetilde{\mathscr{X}}-\mathscr{X}^{(1)})^T
*_M (\mathscr{A}^T *_M \mathscr{R}^{(1)}
      +\mathscr{R}^{(1)} *_M \mathscr{C}^T)\right)\\
&=tr\left((\mathscr{D}-\mathscr{A} *_M \mathscr{X}^{(1)}
      -\mathscr{X}^{(1)}*_N\mathscr{C})^T
 *_M \mathscr{R}^{(1)}\right)\\
&=\|\mathscr{R}^{(1)}\|^2.\\
\end{split}
\end{equation}

\noindent  Assume that \eqref{eqC11} holds for $k=s$, then
\begin{equation*}
\begin{split}
tr\left((\widetilde{\mathscr{X}}-\mathscr{X}^{(s+1)})^T *_M \mathscr{P}^{(s)}\right)
&=tr\left((\widetilde{\mathscr{X}}-\mathscr{X}^{(s)}-\alpha(s) \cdot \mathscr{P}^{(s)})^T
*_M \mathscr{P}^{(s)}\right)\\
&=tr\left((\widetilde{\mathscr{X}}-\mathscr{X}^{(s)})^T *_M \mathscr{P}^{(s)}\right)
-\alpha(s) \cdot tr\left({\mathscr{P}^{(s)}}^T *_M \mathscr{P}^{(s)}\right)\\
&=\|\mathscr{P}^{(s)}\|^2
-\alpha(s) \cdot tr\left({\mathscr{P}^{(s)}}^T *_M \mathscr{P}^{(s)}\right)\\
&=0.\\
\end{split}
\end{equation*}
Furthermore, we have
\begin{equation}\label{eqC13}
\begin{split}
&tr\left((\widetilde{\mathscr{X}}-\mathscr{X}^{(s+1)})^T *_M \mathscr{P}^{(s+1)}\right)\\
&=tr\left((\widetilde{\mathscr{X}}-\mathscr{X}^{(s+1)})^T
*_M
(\mathscr{A}^T*_M \mathscr{R}^{(s+1)}+\mathscr{R}^{(s+1)}*_N \mathscr{C}^T
   +\beta(s+1) \cdot \mathscr{P}^{(s)})\right)\\
&=tr\left((\mathscr{A} *_M (\widetilde{\mathscr{X}}-\mathscr{X}^{(s+1)})
   +(\widetilde{\mathscr{X}}-\mathscr{X}^{(s+1)}) *_N\mathscr{C})^T
*_M \mathscr{R}^{(s+1)}\right)\\
&=tr\left((\mathscr{D} - \mathscr{A} *_M \mathscr{X}^{(s+1)}-\mathscr{X}^{(s+1)} *_N\mathscr{C})^T
*_M \mathscr{R}^{(s+1)}\right)\\
&=\|\mathscr{R}^{(s+1)}\|^2.\\
\end{split}
\end{equation}
The proof is complete. $\square$
\medskip\\
\indent Making use of Lemmas \ref{le3.2} and \ref{le3.3}, we can prove the main result of this paper.
\begin{Theorem}\label{th3.1} If the tensor equation \eqref{eqA02} is consistent, then
for any initial iteration tensor $\mathscr{X}^{(1)}$, its solution can be derived
by Algorithm 3.1 within finite iteration steps.
\end{Theorem}
\noindent\textit{Proof}. For simplicity, denote
$$m:=I_1\cdot \ldots \cdot I_M,\ n:=J_1\cdot \ldots \cdot J_N.$$
If $\mathscr{R}^{(k)}\neq 0$, $k=1, 2, \ldots, mn$, it follows from Lemma \ref{le3.3} that
$\mathscr{P}^{(k)} \neq 0$, then one can compute $\mathscr{X}^{(mn+1)}$
and $\mathscr{R}^{(mn+1)}$ by Algorithm 3.1. Furthermore, from Lemma \ref{le3.2} we know that
$$tr\left({\mathscr{R}^{(mn+1)}}^T *_M \mathscr{R}^{(k)}\right)=0\
\textrm{and}\  tr\left({\mathscr{P}^{(mn+1)}}^T *_M \mathscr{P}^{(l)}\right)=0,$$
where $k,l=1, 2, \ldots, mn,\,k\neq l$.
   Since the sequence $\{\mathscr{R}^{(k)}\}$ is an orthogonal basis of tensor space
$\mathbb{R}^{I_1 \times \dots \times I_M \times J_1 \times \dots \times J_N}$,
which implies that
$\mathscr{R}^{(mn+1)}=0$, i.e., $\mathscr{X}^{(mn+1)}$ is a solution of \eqref{eqA02}. $\square$
\medskip\\
\indent Moreover, according the basic properties of Algorithm 3.1
mentioned above, we can show that
the solvability of the tensor equation \eqref{eqA02} can be determined automatically
during the iteration process.
\begin{Theorem}\label{th3.2}  The tensor equation \eqref{eqA02} is inconsistent
if and only if there exists a positive integer $k_0$ such that $\mathscr{R}^{(k_0)}\neq 0$ and
$\mathscr{P}^{(k_0)}=0$.
\end{Theorem}
\noindent\textbf{Proof}.
If the tensor equation \eqref{eqA02} is inconsistent, it follows that
$\mathscr{R}^{(k)}\neq 0$ for any $k$. Provided that $\mathscr{P}^{(k)}\neq 0$ for all
positive integer $k$, then, from the proof of Theorem \ref{th3.1} we know that
there must exist $\mathscr{X}^{(k)}$ satisfying \eqref{eqA02},
which contradicts to the inconsistency.
   Conversely, if there is a positive integer $k_0$,
such that $R_{k_0}\neq 0$ but $P_{k_0}=0$, which contradicts to Lemma \ref{le3.3},
so the tensor equation \eqref{eqA02} is inconsistent.
   The proof is complete. $\square$

   In addition, since the tensor equation is always over-determined,
we are often interested in the least F-norm solution.
Next we can show that the least F-norm solution of the tensor equation \eqref{eqA02}
can also be gained by means of Algorithm 3.1.
   We first prove the following lemma for this aim.
\begin{Lemma}\label{le3.4}
Let $\mathscr{X}^*$ be a solution of the tensor equation \eqref{eqA01},
then $\mathscr{X}^*$ is the unique least F-norm solution if
$\mathscr{X}^*\in R(\mathscr{A}^T)$.
\end{Lemma}
\noindent\textit{Proof}. For convenience of expression, we use the same symbol $\psi$
to represent the unfoldings of different tensors, e.g.,
$\textbf{A}=\psi(\mathscr{A})$, $\textbf{B}=\psi(\mathscr{B})$ and
$\textbf{X}=\psi(\mathscr{X})$. We prove the conclusion by two steps:

   Step 1) The tensor equation \eqref{eqA01} is equivalent to the matrix equation
\begin{equation}\label{eqC16}
\textbf{AX}=\textbf{B}\ \textrm{with}\ \textbf{X} \in\mathbb{R}^{m\times n}.\\
\end{equation}
   In fact, from the definition of Einstein product, we can respectively
rewrite \eqref{eqA01} and \eqref{eqC16} in components as
$$(\mathscr{A}*_M\mathscr{X})_{k_1 \dots k_P j_1 \dots j_N}
=\sum_{i_1, \dots, i_M}
 \mathscr{A}_{k_1 \dots k_P i_1 \dots i_M}\mathscr{X}_{i_1 \dots i_M j_1 \dots j_N}
=\mathscr{B}_{k_1 \dots k_P j_1 \dots j_N},
$$
and
$$\sum\limits_{t}\textbf{A}_{pt}\textbf{X}_{ts}=\textbf{B}_{ps}.$$
Since $\psi$ is a bijection, then there must exist, respectively,
the unique index
$\{k_1, \ldots, k_P\}$, $\{i_1,\ldots, i_M\}$ and $\{j_1,\ldots, j_N\}$
such that $\textrm{ivec}([k_1, \ldots, k_P],\mathbb{K})=p$,
$\textrm{ivec}([i_1,\ldots, i_M],\mathbb{I})=t$, and
$\textrm{ivec}([j_1,\ldots, j_N],\mathbb{J})=s$ in which $\mathbb{K}=\{K_1,\ldots,K_P\}$.
Therefore, the above two systems are equivalent.

   Step 2) As is well-known [21], the least F-norm solution of matrix equation \eqref{eqC16} is
$\widetilde{\textbf{X}}=\textbf{A}^\dagger \textbf{B}\in R(\textbf{A}^T)$,
where the superscript $\dagger$ denotes the Moore-Penrose inverse of a matrix.
   In view of the uniqueness of the least F-norm solution,
and together with the fact that $\psi$ is a bijection,
we complete the proof. $\square$

  Depending on the above lemma, we can prove the following theorem.

\begin{Theorem}\label{th3.3} Assume that the tensor equation \eqref{eqA02} is consistent,
and let the initial iteration tensor
$\mathscr{X}^{(1)}=\mathscr{A}^T *_M \mathscr{W}+\mathscr{W} *_N \mathscr{C}^T$ with arbitrary
$\mathscr{W}\in \mathbb{R}^{I_1\times \dots \times I_M \times J_1\times \dots \times J_N}$,
or especially, $\mathscr{X}^{(1)}=\mathscr{O}
\in \mathbb{R}^{I_1\times \dots \times I_M \times J_1\times \dots \times J_N}$,
then the solution generated by Algorithm 3.1 is the unique least F-norm solution.
\end{Theorem}
\medskip
\noindent\textit{Proof}. From Algorithm 3.1 and Theorem 3.1, it is known
that if we choose the initial iteration tensor
$\mathscr{X}^{(1)}=\mathscr{A}^T *_M \mathscr{W}+\mathscr{W} *_N \mathscr{C}^T$
for some tensor $\mathscr{W}$, then the approximate solution $\mathscr{X}^{(k)}$
of \eqref{eqA02} possesses the form
$\mathscr{X}^{(k)}=\mathscr{A}^T *_M \mathscr{H} + \mathscr{H} *_N \mathscr{C}^T$
for some
$\mathscr{H}\in \mathbb{R}^{I_1\times \dots \times I_M \times J_1\times \dots \times J_N}$.
   Using the definition of Vec and Lemma \ref{le2.1}, we deduce that
\begin{equation}\label{eqC14}
\begin{split}
\textrm{Vec}\left(\mathscr{X}^{(k)}\right)
&=\textrm{Vec}\left(\mathscr{A}^T *_M\mathscr{H}+\mathscr{H} *_N \mathscr{C}^T\right)\\
&=\left(\mathscr{I}_1\otimes\mathscr{A}^T+\mathscr{C}\otimes\mathscr{I}_2\right) *_N \textrm{Vec}\left(\mathscr{H}\right)\\
&\in R\left((\mathscr{I}_1\otimes\mathscr{A}+\mathscr{C}^T\otimes\mathscr{I}_2)^T\right),\\
\end{split}
\end{equation}
where $\mathscr{I}_1$ and $\mathscr{I}_2$ are the identity tensors of
size $J_1\times \dots \times J_N \times J_1\times \dots \times J_N$ and
$I_1\times \dots \times I_M \times I_1\times \dots \times I_M$, respectively.
   On the other hand, by using the properties of the Kronecker product,
one can demonstrate that \eqref{eqA02} is equivalent to the tensor equation
\begin{equation*}\label{eqC15}
\begin{split}
(\mathscr{I}_1\otimes\mathscr{A}+\mathscr{C}^T\otimes\mathscr{I}_2) *_N \textrm{Vec}(\mathscr{X})=\textrm{Vec}(\mathscr{D}),\\
\end{split}
\end{equation*}
which, together with \eqref{eqC14} and Lemma \ref{le3.4}, implies that
$\mathscr{X}^{(k)}$ is the least F-norm solution of the Sylvester tensor equation \eqref{eqA02}.
The proof is complete. $\square$

\section{Solving the tensor nearness problem}
   In this section, we apply Algorithm 3.1 to the solution
of the tensor nearness problem \eqref{eqA03}.
   Suppose that the tensor equation \eqref{eqA02} is consistent,
i.e., its solution set, denoted by $\Phi$, is nonempty.
   It is easy to verify that the set $\Phi$ is a closed and convex set in the tensor space
$\mathbb{R}^{I_1\times \dots \times I_M \times J_1\times \dots \times J_N}$,
which reveals that the solution to the tensor nearness problem is unique,
denoted by $\widehat{\mathscr{X}}$ for convenience.

   We should point out that the unique solution $\widehat{\mathscr{X}}$ can also be
derived by using Algorithm 3.1.
   Actually, noting the fact that
to solve the tensor nearness problem with the given tensor $\mathscr{X}_0$ is equivalent to
find the least F-norm solution (denoted by $\widehat{\mathscr{Y}}$)
of the following Sylvester tensor equation
\begin{equation}\label{eqD01}
\mathscr{A}*_M\mathscr{Y}+\mathscr{Y}*_N\mathscr{C}=\widetilde{\mathscr{D}},\\
\end{equation}
where $\mathscr{Y}=\mathscr{X}-\mathscr{X}_0$ and
$\widetilde{\mathscr{D}}=\mathscr{D}-\mathscr{A}*_M\mathscr{X}_0-\mathscr{X}_0*_N\mathscr{C}$,
then, it follows from Theorem \ref{th3.3} that $\widehat{\mathscr{X}}$
can be obtained by applying Algorithm 3.1 to \eqref{eqD01} with the initial iteration tensor
$\mathscr{X}^{(1)}=\mathscr{A}^T *_M \mathscr{W}+\mathscr{W} *_N \mathscr{C}^T$ for some
$\mathscr{W}\in \mathbb{R}^{I_1\times \dots \times I_M \times J_1\times \dots \times J_N}$,
or especially, $\mathscr{X}^{(1)}=\mathscr{O}
\in \mathbb{R}^{I_1\times \dots \times I_M \times J_1\times \dots \times J_N}$.
In this case, the nearness solution of \eqref{eqA03} can be obtained by
$\widehat{\mathscr{X}}=\widehat{\mathscr{Y}}+\mathscr{X}_0$.

\section{Numerical experiments}

In this section, we perform some numerical examples to illustrate the feasibility and
effectiveness of the proposed algorithm in present paper.
All computations were written using MATLAB (version R2016a) on a personal computer with 2.50GHz
central processing unit (Intel(R) Core(TM) i5-3210M) and 4GB memory.
Specially, all the tensor calculations in our tests were carried out with the Tensor Toolbox Version
2.6.\footnote{http://www.sandia.gov/tgkolda/TensorToolbox/index-2.6.html.}
   The iterations will be terminated if the norm of the residual, i.e.,
RES$=\|\mathscr{D}-\mathscr{A}*_M\mathscr{X}^{(k)}-\mathscr{X}^{(k)}*_N\mathscr{C}\|
<\varepsilon=1.0e-10$, or
the number of iteration steps exceeds the maximum $k_{\max}=1000$.

\begin{Example}\label{ex01}
Suppose the tensors $\mathscr{A}\in \mathbb{R}^{4\times 3\times 4\times 3}$,
$\mathscr{C}\in \mathbb{R}^{3\times 3\times 3\times 3}$,
$\mathscr{D}\in \mathbb{R}^{4\times 3\times 3\times 3}$ in \eqref{eqA02}
are given as follows:
$$
\mathscr{A}(:,:,1,1)=
\left[
  \begin{array}{ccc}
	    11&     7&     7\\
	    -2&    11&    -2\\
	    11&    -2&     7\\
	    -2&    11&    -2\\
  \end{array}
\right], \mathscr{A}(:,:,2,1)=
\left[
  \begin{array}{ccc}
	    -2 &   -2&    -2\\
	     3 &   -2&     3\\
	    -2 &    3&    -2\\
	     3 &   -2&     3\\
  \end{array}
\right],
$$
$$
\mathscr{A}(:,:,3,1)=
\left[
  \begin{array}{ccc}
	     3 &   -4 &   -4\\
	    -1 &    3 &   -1\\
	     3 &   -1 &   -4\\
	    -1 &    3 &   -1\\
  \end{array}
\right],
\mathscr{A}(:,:,4,1)=
\left[
  \begin{array}{ccc}
	     2 &   -9 &   -9\\
	    -6 &    2 &   -6\\
	     2 &   -6 &   -9\\
	    -6 &    2 &   -6\\
  \end{array}
\right],
$$
$$
\ \mathscr{A}(:,:,2,2)=
\left[
  \begin{array}{ccc}
	   -16 &    3 &    3\\
	   -11 &  -16 &  -11\\
	   -16 &  -11 &    3\\
	   -11 &  -16 &  -11\\
  \end{array}
\right],
\mathscr{A}(:,:,1,2)=
\left[
  \begin{array}{ccc}
	     0 &    7 &    7\\
	    11 &    0 &   11\\
	     0 &   11 &    7\\
	    11 &    0 &   11\\
  \end{array}
\right],
$$
$$\
\mathscr{A}(:,:,3,2)=
\left[
  \begin{array}{ccc}
	   -11&    15 &   15\\
	     0&   -11 &    0\\
	   -11&     0 &   15\\
	     0&   -11 &    0\\
  \end{array}
\right],
\mathscr{A}(:,:,4,2)=
\left[
  \begin{array}{ccc}
	    -4 &   -2 &   -2\\
	    16 &   -4 &   16\\
	    -4 &   16 &   -2\\
	    16 &   -4 &   16\\
  \end{array}
\right],
$$
$$\
\mathscr{A}(:,:,1,3)=
 \left[
   \begin{array}{ccc}
	     3 &   -3 &   -3\\
	    13 &    3 &   13\\
	     3 &   13 &   -3\\
	    13 &    3 &   13\\
   \end{array}
 \right],
 \mathscr{A}(:,:,2,3)=
 \left[
   \begin{array}{ccc}
	    26&     0 &    0\\
	    -4&    26 &   -4\\
	    26&    -4 &    0\\
	    -4&    26 &   -4\\
   \end{array}
 \right],\ \ \
$$
$$
\mathscr{A}(:,:,3,3)=
 \left[
   \begin{array}{ccc}
	    -4 &    1&     1\\
	     8 &   -4&     8\\
	    -4 &    8&     1\\
	     8 &   -4&     8\\
   \end{array}
 \right],
 \mathscr{A}(:,:,4,3)=
 \left[
   \begin{array}{ccc}
	     2  &  -8&    -8\\
	   -16  &   2&   -16\\
	     2  & -16&    -8\\
	   -16  &   2&   -16\\
   \end{array}
 \right];
$$

$$\ \
\mathscr{C}(:,:,1,1)=
\left[
  \begin{array}{ccc}
	    10  &   0 &    6\\
	    15  &  10 &   10\\
	    10  &  15 &   10\\
\end{array}
\right], \
\mathscr{C}(:,:,2,1)=
\left[
  \begin{array}{ccc}
	     6 &   -9 &   17\\
	    -9 &    6 &    6\\
	     6 &   -9 &    6\\
\end{array}
\right], \ \ \ \ \ \
$$
$$
\mathscr{C}(:,:,3,1)=
\left[
  \begin{array}{ccc}
	     4 &  -19 &   -3\\
	   -14 &    4 &    4\\
	     4 &  -14 &    4\\
\end{array}
\right],  \ \
\mathscr{C}(:,:,1,2)=
\left[
  \begin{array}{ccc}
	     9 &  -22 &   -8\\
	     0 &    9 &    9\\
	     9 &    0 &    9\\
\end{array}
\right],
$$
$$\
\mathscr{C}(:,:,2,2)=
\left[
  \begin{array}{ccc}
	     0  &  -9 &   -3\\
	   -13  &   0 &    0\\
	     0  & -13 &    0\\
\end{array}
\right],\
\mathscr{C}(:,:,3,2)=
\left[
  \begin{array}{ccc}
	    -7 &  -17 &   12\\
	     6 &   -7 &   -7\\
	    -7 &    6 &   -7\\
\end{array}
\right],
$$
$$
\mathscr{C}(:,:,2,3)=
\left[
  \begin{array}{ccc}
	     5  & -13 &    1\\
	    -5  &   5 &    5\\
	     5  &  -5 &    5\\
\end{array}
\right],\ \ \ \
 \mathscr{C}(:,:,1,3)=
\left[
  \begin{array}{ccc}
	     0  &  -3 &    4\\
	     5  &   0 &    0\\
	     0  &   5 &    0\\
\end{array}
\right], \ \ \ \ \ \
$$
$$
 \mathscr{C}(:,:,3,3)=
\left[
  \begin{array}{ccc}
	     0  & -12  &   3\\
	    -1  &   0  &   0\\
	     0  &  -1  &   0\\
\end{array}
\right],
$$
and the tensor $\mathscr{D}$ is chosen such that
$\mathscr{D}=\mathscr{A}*_M\mathscr{X}^*+\mathscr{X}^* *_N\mathscr{C}$
with $\mathscr{X}^*=\verb"reshape"(1:108,[4,3,3,3])
\in \mathbb{R}^{4\times 3\times 3\times 3}$.
\end{Example}

\begin{figure}[h]
\centering
\includegraphics[totalheight=4.00in]{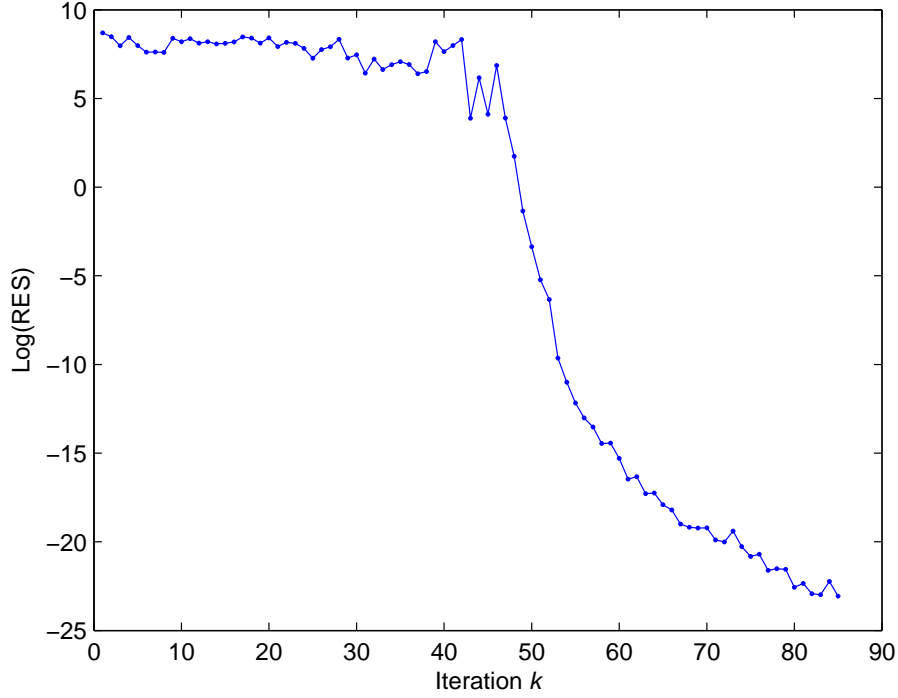}
\centering
\caption{\footnotesize Convergence behavior of Algorithm 3.1 for the tensor equation from Example \ref{ex01}.}\label{fig1}
\end{figure}

   In this case, the tensor equation \eqref{eqA02} is consistent and
$\mathscr{X}^*$ is an exact solution.
   Applying Algorithm 3.1 with initial iteration tensor $\mathscr{X}^{(1)}=\mathscr{O}$
to \eqref{eqA02},
we obtain the least F-norm solution, denoted by $\widetilde{\mathscr{X}}$, and the
corresponding residual RES$=9.6392\textrm{e}-11$
after 86 iteration steps.
   Moreover, the convergence behavior of our algorithm is plotted in Figure \ref{fig1},
which appears that this method is efficient.
   Many other tests not reported here also confirm this phenomenon.
$$
\small{
 \widetilde{\mathscr{X}}(:,:,1,1)=
\left[
  \begin{array}{ccc}
	   42.9784 &  46.3496&   53.2346\\
	   53.0555 &  68.2438&   49.0433\\
	   54.4996 &  54.9506&   59.0609\\
	   61.1020 &  53.8303&   73.3694\\
  \end{array}
\right],
 \widetilde{\mathscr{X}}(:,:,2,1)=
\left[
  \begin{array}{ccc}
	   32.9897&   36.6903 &  42.0641\\
	   38.3122&   47.6399 &  40.5920\\
	   39.5236&   41.8336 &  45.8862\\
	   43.1914&   41.8239 &  53.2235\\
  \end{array}
\right], }
$$
$$
\small{
 \widetilde{\mathscr{X}}(:,:,3,1)=
\left[
  \begin{array}{ccc}
	   46.9887 &  50.6593  & 56.1705\\
	   52.7434 &  62.6039  & 54.4513\\
	   53.9760 &  56.1170  & 60.1748\\
	   57.9106 &  56.0063  & 68.1459\\
  \end{array}
\right],
 \widetilde{\mathscr{X}}(:,:,1,2)=
\left[
  \begin{array}{ccc}
	   37.0000 &  41.0000  & 45.0000\\
	   38.0000 &  42.0000  & 46.0000\\
	   39.0000 &  43.0000  & 47.0000\\
	   40.0000 &  44.0000  & 48.0000\\
  \end{array}
\right], }
$$
$$
\small{
 \widetilde{\mathscr{X}}(:,:,2,2)=
\left[
  \begin{array}{ccc}
	   50.9990 &  54.9690  & 59.1064\\
	   52.4312 &  56.9640  & 59.8592\\
	   53.4524 &  57.2834  & 61.2886\\
	   54.7191 &  58.1824  & 62.9224\\
  \end{array}
\right],
 \widetilde{\mathscr{X}}(:,:,3,2)=
\left[
  \begin{array}{ccc}
	   41.0103  & 45.3097 &  47.9359\\
	   37.6878  & 36.3601 &  51.4080\\
	   38.4764  & 44.1664 &  48.1138\\
	   36.8086  & 46.1761 &  42.7765\\
  \end{array}
\right],}
$$
$$
\small{
 \widetilde{\mathscr{X}}(:,:,1,3)=
 \left[
   \begin{array}{ccc}
	   73.0000  & 77.0000 &  81.0000\\
	   74.0000  & 78.0000 &  82.0000\\
	   75.0000  & 79.0000 &  83.0000\\
	   76.0000  & 80.0000 &  84.0000\\
   \end{array}
 \right],
 \widetilde{\mathscr{X}}(:,:,2,3)=
 \left[
   \begin{array}{ccc}
	   57.0144 &  61.4336 &  63.5103\\
	   51.9630 &  48.5042 &  67.9711\\
	   52.6669 &  59.0329 &  62.9594\\
	   49.9320 &  61.4465 &  55.0871\\
   \end{array}
 \right], }
$$
$$
\small{
 \widetilde{\mathscr{X}}(:,:,3,3)=
 \left[
   \begin{array}{ccc}
	   59.0196 &  63.5885 &  64.9782\\
	   51.8069 &  45.6842 &  70.6751\\
	   52.4051 &  59.6161 &  63.5163\\
	   48.3363 &  62.5345 &  52.4753\\
   \end{array}
 \right].}
$$

\newpage
   Next we consider the tensor nearness problem \eqref{eqA03}.
\begin{Example}\label{ex02}
Let the tensors $\mathscr{A}$, $\mathscr{C}$ and $\mathscr{D}$ in \eqref{eqA02}
be the same as in Example \ref{ex01}, and assume that the given tensor
$\mathscr{X}_0$ is as follows:
$$
\small{
 \mathscr{X}_0(:,:,1,1)=
\left[
  \begin{array}{ccc}
	     0 &   11 &   -10\\
	    -7 &   -1 &    -4\\
	    -4 &   -4 &     5\\
	     7 &   -6 &    -5\\
  \end{array}
\right],
 \mathscr{X}_0(:,:,2,1)=
\left[
  \begin{array}{ccc}
	     4 &   -11 &   -12\\
	    -3 &     6 &   -20\\
	   -13 &     4 &     0\\
	     6 &    -6 &    -4\\
  \end{array}
\right], }
$$
$$
\small{
 \mathscr{X}_0(:,:,3,1)=
\left[
  \begin{array}{ccc}
	    -5  &   11 &     0\\
	   -16  &   -2 &     4\\
	    33  &   -2 &    -1\\
	   -16  &   -8 &     9\\
  \end{array}
\right],
 \mathscr{X}_0(:,:,1,2)=
\left[
  \begin{array}{ccc}
	     7  &    6&     -4\\
	    14  &   -4&    -11\\
	   -13  &  -32&      9\\
	   -28  &    0&    -10\\
  \end{array}
\right], }
$$
$$
\small{
 \mathscr{X}_0(:,:,2,2)=
\left[
  \begin{array}{ccc}
	   -10  &    5 &    18\\
	    -6  &   -8 &     8\\
	   -16  &   -4 &     8\\
	   -12  &   -4 &     9\\
  \end{array}
\right], \
 \mathscr{X}_0(:,:,3,2)=
\left[
  \begin{array}{ccc}
	    -7  &   11 &     4\\
	    -4  &   -1 &     0\\
	   -14  &   -5 &   -21\\
	     4  &    6 &    14\\
  \end{array}
\right],\ }
$$
$$
\small{
 \mathscr{X}_0(:,:,1,3)=
 \left[
   \begin{array}{ccc}
	     1  &    1 &     4\\
	     7  &   21 &    -5\\
	    -4  &    2 &     0\\
	   -16  &    5 &   -18\\
   \end{array}
 \right],
\mathscr{X}_0(:,:,2,3)=
 \left[
   \begin{array}{ccc}
	     9  &    4  &    5\\
	     2  &   -2  &   -4\\
	    -7  &   -2  &   13\\
	   -16  &    6  &   -4\\
   \end{array}
 \right],\ \  }
$$
$$
\small{
\mathscr{X}_0(:,:,3,3)=
 \left[
   \begin{array}{ccc}
	     0 &    -9 &     1\\
	     8 &   -16 &   -14\\
	     9 &   -15 &   -12\\
	   -19 &    -3 &    -2\\
   \end{array}
 \right].}
$$
\end{Example}

   Applying Algorithm 3.1 with $\mathscr{X}^{(1)}=\mathscr{O}$ to the tensor equation \eqref{eqD01},
we obtain the solution to the tensor nearness problem \eqref{eqA03}
after 79 iteration steps,
i.e., $\widehat{\mathscr{X}}$.
$$
\small{
 \widehat{\mathscr{X}}(:,:,1,1)=
\left[
  \begin{array}{ccc}
	   44.1912&   54.1943&   43.7075\\
	   50.4602&   68.6229&   49.3393\\
	   48.1731&   53.1240&   62.4182\\
	   79.5249&   53.7313&   65.1154\\
  \end{array}
\right],
 \widehat{\mathscr{X}}(:,:,2,1)=
\left[
  \begin{array}{ccc}
	   39.4807&   25.7036&   37.5108\\
	   39.1969&   44.5352&   36.9080\\
	   41.6060&   39.9043&   56.2875\\
	   45.9264&   37.6248&   44.2264\\
  \end{array}
\right], }
$$
$$
\small{
 \widehat{\mathscr{X}}(:,:,3,1)=
\left[
  \begin{array}{ccc}
	   40.5057&   59.6836&   56.9232\\
	   41.8432&   65.1252&   55.9031\\
	   83.6837&   59.5245&   59.2572\\
	   52.0225&   56.2288&   71.7468\\
  \end{array}
\right],
 \widehat{\mathscr{X}}(:,:,1,2)=
\left[
  \begin{array}{ccc}
	   37.0000 &  41.0000  & 45.0000\\
	   38.0000 &  42.0000  & 46.0000\\
	   39.0000 &  43.0000  & 47.0000\\
	   40.0000 &  44.0000  & 48.0000\\
  \end{array}
\right], }
$$
$$
\small{
 \widehat{\mathscr{X}}(:,:,2,2)=
\left[
  \begin{array}{ccc}
	   41.7494&   54.6924&   71.7399\\
	   49.8603&   58.7485&   67.0232\\
	   33.7992&   61.6395&   66.1633\\
	   53.6123&   55.2035&   71.7550\\
  \end{array}
\right],
 \widehat{\mathscr{X}}(:,:,3,2)=
\left[
  \begin{array}{ccc}
	   34.5193 &  56.2964&   52.4892\\
	   36.8031 &  39.4648&   55.0920\\
	   36.3940 &  46.0957&   37.7125\\
	   34.0736 &  50.3752&   51.7736\\
  \end{array}
\right],}
$$
$$
\small{
 \widehat{\mathscr{X}}(:,:,1,3)=
 \left[
   \begin{array}{ccc}
	   73.0000  & 77.0000 &  81.0000\\
	   74.0000  & 78.0000 &  82.0000\\
	   75.0000  & 79.0000 &  83.0000\\
	   76.0000  & 80.0000 &  84.0000\\
   \end{array}
 \right],
 \widehat{\mathscr{X}}(:,:,2,3)=
 \left[
   \begin{array}{ccc}
	   66.1178 &  55.6214 &  63.4721\\
	   57.5053 &  53.0468 &  68.9621\\
	   44.5545 &  56.3846 &  70.5580\\
	   50.0412 &  66.5350 &  56.9444\\
   \end{array}
 \right], }
$$
$$
\small{
 \widehat{\mathscr{X}}(:,:,3,3)=
 \left[
   \begin{array}{ccc}
	   64.4359 &  52.8083 &  61.1574\\
	   62.3311 &  36.4565 &  60.7722\\
	   56.7895 &  56.3274 &  48.6033\\
	   36.7992 &  60.3013 &  46.4383\\
   \end{array}
 \right].}
$$
At this time, $\|\widehat{\mathscr{X}}-\mathscr{X}_0\|=640.2422$.
\section{Conclusions}

   In this paper, we present an iterative method for solving the Sylvester tensor equation \eqref{eqA02},
i.e., Algorithm 3.1.
   For any initial iteration tensor, it is shown that the solvability of this equation
can be determined automatically (see, Theorem \ref{th3.2}),
and that the solution (if it exits) can be obtained within
finite iteration steps in absence of roundoff errors (see, Theorem \ref{th3.1}).
   Particularly, the least F-norm solution of \eqref{eqA02} can also be derived by selecting
appropriate initial iteration tensor (see, Theorem \ref{th3.3}).
   Additionally, applying this iterative method to another Sylvester tensor equation, i.e., \eqref{eqD01},
we can obtain the unique solution to the tensor nearness problem \eqref{eqA03}.
   Many other examples we have tested in MATLAB confirm the
theoretical results presented in this paper.
   Of course, for a problem with large and not sparse
tensors $\mathscr{A}, \mathscr{C}$ and $\mathscr{D}$,
Algorithm 3.1 may not terminate in a finite number of iteration steps because
of roundoff errors. This is an important problem which we should study in a future work.
   Moreover, the approach we propose in this paper can not be directly used to solve
the Sylvester tensor equation \eqref{eqA02} when it is inconsistent,
which will be considered in our future work as well.

%

\end{document}